\documentclass{amsart}

\usepackage{mathrsfs}
\usepackage{multicol}
\usepackage{pdfpages}
\usepackage{wrapfig}
\usepackage{color}
\usepackage{bbm}
\usepackage{tabularx}
\usepackage{caption}
\usepackage{bm}
\usepackage{mathrsfs}
\usepackage{amssymb}
\usepackage{enumitem}
\usepackage{pdfpages}
\usepackage{float}
\usepackage{graphicx}
\usepackage{csquotes}
\usepackage{hyperref}
\usepackage{amsmath}
\usepackage{ulem}

\newtheorem{theo}{Theorem}

\newtheorem{lem}{Lemma}

\usepackage{fancyhdr}
\usepackage{ulem}
\usepackage{mathrsfs}
\usepackage[english]{babel}
\usepackage{multicol}
\usepackage{pdfpages}
\usepackage{wrapfig}
\usepackage{color}
\usepackage{tabularx}
\usepackage{caption}
\usepackage{bm}
\usepackage{mathrsfs}
\usepackage{amssymb}
\usepackage{enumitem}
\usepackage{pdfpages}
\usepackage{float}
\usepackage{graphicx}
\usepackage{csquotes}
\usepackage{hyperref}
\usepackage{amsmath}
\usepackage{ulem}
\usepackage{lipsum}

\theoremstyle{definition}

\theoremstyle{remark}

\numberwithin{equation}{section}

\numberwithin{equation}{section}
\numberwithin{theorem}{section}
\numberwithin{cor}{section}
\numberwithin{cor}{section}
\numberwithin{eg}{section}
\numberwithin{examp}{section}

\newcommand{\A}{{\mathcal A}}
\newcommand{\bbr}{{\mathcal R}}
\newcommand{\bbn}{{\mathbb N}}
\newcommand{\raw}{\rightarrow}

\newcommand{\nti}{n\rightarrow\infty}

\begin{document}

\title{Central Limit Theorem in High Dimensions : The Optimal Bound on Dimension Growth Rate}

%\title{\rd{The Critical growth rate of dimension in high dimensional Central Limit Theorem}}

\author{Debraj Das}
\address{Department of Mathematics and Statistics, Indian Institute of Technology, Kanpur, Uttar Pradesh 208016, India}
%\curraddr{Department of Mathematics and Statistics, Washington University in St. Louis, MO 63130, United States}
\email{rajdas@iitk.ac.in}
\thanks{The first author was supported in part by DST fellowship DST/INSPIRE/04/2018/001290}

\author{Soumendra Lahiri}
\address{Department of Mathematics and Statistics, Washington University in St. Louis, MO 63130, United States}
\email{s.lahiri@wustl.edu}
\thanks{The second author was supported in part by NSF grant DMS 2006475}

%\authorrunning{Das \and Lahiri}
%\lhead{Owl stretching time}

%    General info

%\date{January 1, 2001 and, in revised form, June 22, 2001.}



\begin{abstract}
In this article, we try to give an answer to the simple question: ``\textit{What is the critical growth rate of the dimension $p$ as a function of the sample size $n$ for which the Central Limit Theorem holds uniformly 
over the collection of $p$-dimensional hyper-rectangles ?''}. 
Specifically, we are interested in the normal approximation of suitably scaled versions  of the sum $\sum_{i=1}^{n}X_i$ in $\bbr^p$ uniformly over the class of hyper-rectangles $\mathcal{A}^{re}=\{\prod_{j=1}^{p}[a_j,b_j]\cap\mathcal{R}:-\infty\leq a_j\leq b_j
\leq \infty, j=1,\ldots,p\}$, where $X_1,\dots,X_n$ are independent $p-$dimensional random vectors with each having independent and identically distributed (iid) components. We investigate the critical  cut-off rate of $\log p$ below which the uniform central limit theorem (CLT) holds and above which it fails. According to some recent results of Chernozukov et al. (2017), it is well known that the CLT holds uniformly over $\mathcal{A}^{re}$ if $\log p=o\big(n^{1/7}\big)$. They also conjectured that for CLT to hold uniformly over $\mathcal{A}^{re}$, the optimal rate is $\log p = o\big(n^{1/3}\big)$. We show instead that under some conditions, the CLT holds uniformly over $\mathcal{A}^{re}$, when $\log p=o\big(n^{1/2}\big)$. More precisely, we show that  if  $\log p =\epsilon \sqrt{n}$ for some sufficiently small $\epsilon>0$, the normal approximation is valid with an error $\epsilon$, uniformly over $\mathcal{A}^{re}$. Further, we show by an example that the uniform CLT over $\mathcal{A}^{re}$ fails if $\limsup_{\nti}
n^{-(1/2+\delta)}
\log p >0$ for some $\delta>0$.
Hence the critical  rate of the growth of $p$ for the validity of the 
CLT is given by $\log p=o\big(n^{1/2}\big)$.
%\rd{something in between $\sqrt{n}$ and $\sqrt{n\log n}$.}
%Indeed we solved and extended the %conjecture of Chernozukov et al. (2017) %for the class of hyper-rectangles. 
\end{abstract}

\maketitle

\section{Introduction}\label{sec:intro}
{Central Limit Theorem (CLT) is one of the oldest as well as remarkable results of classical probability theory. After initial works by de Moivre in the eighteenth century and by Laplace in the nineteenth century, it is the second half of the twentieth century which sees a boom in different forms as well as different applications of the CLT. In most simplest words, CLT is a statement about the convergence of properly centered and scaled sample mean of a sequence of random vectors to the Gaussian random vector in distribution. Although most of the theoretical developments centred around to establish CLT whenever the underlying dimension of the random vector is fixed, the recent interest, primarily driven by problems arising in statistical inference and machine learning,
lies in establishing CLT when the dimension also grows with the sample size. Hence a natural but important question is ``\textit{What is the critical growth rate of dimension $p$ as a function of the sample size $n$ for the validity of the CLT
in the high dimensional set up where $p\gg n$?}'' Let us consider a simple example to motivate the problem. Suppose that there are a collection of $np$ independent and identically distributed (iid) Rademacher random variables, all defined on the same probability space. Now assume that we arrange those $np$ random variables in $n$ many vectors each of length $p$. Let us denote those $n$ random vectors as $\{{Y}_1,\dots,{Y}_n\}$ and let ${W}_n=n^{-1/2}\sum_{i=1}^{n}{Y}_i$. Now if $p$ is fixed, then classical CLT implies that ${W}_n$ converges in distribution to $Z$ where $Z$ is random vector with the Gaussian distribution with mean zero and covariance matrix ${\mathbb I}_p$, the identity matrix of order $p$. When $p$ grows with $n$, the question is how large $p$ can be compared to $n$ for the Gaussian approximation to hold and obviously what is the critical growth rate of $p$ above which the Gaussian approximation fails. This paper centers around these two questions, but obviously under a more general framework which  we will describe  below.}

Let $X_1,\dots,X_n$ be  independent 
 random vectors in $\mathcal{R}^p$, 
$p\in {\mathbb N}\equiv 
 \{1,2,\ldots\}$ 
 and  let $S_n=X_1+\ldots+X_n$, $n\in {\mathbb N}$.
 Suppose that $EX_i=0$  and $E\|X_i\|^2<\infty$
 for all $i$ where $\|\cdot\|$ denotes the Euclidean norm
 on $\bbr^p$. 
 Lindeberg's  Central Limit Theorem (CLT) (cf. Theorem 11.1.1, Athreya and Lahiri (2006)) asserts that for $p$ fixed (i.e., not changing with $n$),
 under a mild condition on the truncated 
 second moments, 
 \begin{equation}
 \label{clt-1}
  T_n \equiv \Big(\sum_{i=1}^n EX_iX_i'\Big) ^{-1/2}S_n
 \Rightarrow Z
 \end{equation}
 where $\Rightarrow$ denotes 
 convergence in distribution and $Z$ is random vector with the Gaussian distribution 
with mean zero and covariance matrix
${\mathbb I}_p$,
the identity matrix of order $p$. Here and in the following, $B^{\prime}$ denotes the transpose of a matrix $B$. This 
yields the Gaussian approximation 
\begin{equation}
  \rho_{n,\A}\equiv   \sup_{A\in \A} \Big| P(T_n\in A) - P(Z\in A)   \Big| \raw 0  ~\mbox{as}~ \nti,
    \label{clt-2}
\end{equation}
where $\A$ is a suitable collection of convex sets  
in $\bbr^p$. Typical choices of $\A$  include 
\begin{enumerate}[label=(\roman*)]
\item
$\A^{dist}
= \Big\{ (-\infty, a_1]\times\ldots\times (-\infty, a_p]: a_1,\ldots, a_p\in \bbr\Big\}$,\\
the collection of all left-infinite rectangles, leading to the Kolmogorov distance between the distributions of $T_n$ and $Z$, 
\item 
$\A^{max} =\Big\{
 (-\infty, t]\times\ldots\times (-\infty, t]: t \in \bbr
\Big\}
 = \Big\{ \{max_{1\leq j\leq p} T_{nj} \leq t\} : t\in \bbr
 \Big\}$,
~~and
\item
$\mathcal{A}^{re}=\Big\{\prod_{j=1}^{p}[a_j,b_j]\cap\mathcal{R}:-\infty\leq a_j\leq  b_j\leq \infty
~\mbox{for}~ j=1,\ldots,p \Big\}$,\\
the collection of 
all hyper rectangles, 
\end{enumerate}
among others. Clearly, $\A^{max}\subset \A^{dist}\subset \A^{re}$. 
For a fixed $p\in \bbn$, \eqref{clt-1} implies that 
 $\rho_{n,\A^{re}} \raw 0$, so that the Gaussian approximation (GA) holds for each of the three classes. 
In recent years there has been a surge of  interest in 
extending the GA results to the case
 where $p=p_n\rightarrow \infty$ as 
 $n\rightarrow\infty$. In this paper, we investigate the 
 range of validity of the approximation \eqref{clt-2}
 for $\A = \A^{re}$ for increasing $p$ under some 
 suitable regularity conditions.

 To provide some perspective on the problem in relation to existing
 work on the CLT, we point out
% (cf.  Chernozukov et al. (2013)) 
 that when $p\raw \infty$ with $n$,  the 
 class of sets allowed in $\A^{\max}$ (and hence in $\A^{dist}$
 and $\A^{re}$) in our framework does  not necessarily allow the GA problem to
 be embedded in the paths of an empirical process or in some suitable Banach space
 and hence can not be  directly derived from  the well developed body of work
 establishing the CLT therein
 (cf. Ledoux and Talagrand (1991) and van der Vaart and Wellner (2000)). 
 As a result, alternative theoretical tools  are needed. We first highlight some related work and associated approaches that have been applied to study 
 the GA problem in our framework.  
 Portnoy (1986)
 obtained some early results
 on the CLT in increasing dimensions  using 
 Fourier transform techniques allowing
 $p$ to increase as a fractional power of
 $n$. Similar growth rates of $p$ were
 also  allowed in the works of Nagaev (1976), 
 Asriev and Rotar  (1989), and 
 G\"otze (1991)  in 
 studying Berry-Esseen type bounds on the rate of convergence in the CLT
 for 
 different classes of sets and functions 
 under varying degree of 
 generality.  Tiro (1991) derived Edgeworth expansion results
 for expectations of functions of $T_n$,
 again using Fourier transformation based techniques, but under 
   stronger conditions on $p$.

  Two other widely-used approaches 
  for proving the CLT in high dimensions are based on Stein's method (Stein (1986))
  and on Lindeberg's method (Lindeberg (1922)).   CLTs in high dimensions using 
  Stein's method  have been proved by Goldstein and Rinott  (1996)
  by applying size bias couplings, by Chatterjee and Meckes (2008)
  and Reinert and R\"{o}llin (2009) using exchangeable pairs, and by 
  Chen and Fang (2011) using the concentration inequality approach, among others. 
  Lindeberg's method was revived by Trotter (1959) in the context of
  proving multivariate CLTs and  has been significantly  generalized by 
  Chatterjee (2006) for approximating expectations of  smooth functions
  (not necessarily based on sums) of random vectors in  high dimensions. Some variants of the argument
  also have been used by Bentkus (2003) 
  and Zhilova (2019) to establish Berry-Esseen Theorems 
   in high dimensions for the class of all 
   convex sets
  and for the class of all Euclidean balls, respectively. 
  Building on %inspired by  
  Chatterjee's work and using 
  techniques from spin-glass theory, 
  in a seminal work,  Chernozukov et al. (2013) (hereafter referred to as [CCK]) establish   CLT
  for the class of sets $\A^{max}$ 
  in ultra-high dimensions. More precisely, they showed that
  %under 
  %some very  general regularity conditions, 
  $\rho_{n,\A^{max}}\raw 0$ 
  allowing the dimension $p$ to grow
  at a sub-exponential rate: 
  \begin{equation}
  \rho_{n,\A^{max}}\raw 0, ~\mbox{provided}~ \log p = o(n^{1/7}) ~\mbox{as}~ \nti. 
   \label{CCK-1} 
  \end{equation}
  Extensions of  the GA results to the class $\A^{re}$ 
  has been later proved in 
  Chernozukov et al. (2017). (Hereafter, we shall refer to
  both these papers as [CCK]).
  [CCK] also conjectured that the best  growth rate of $p$ is
  possibly faster, and 
  hypothesized the best rate as $\log p = o(n^{1/3})$. 
  A series of recent papers tried to settle this issue, 
  and extend it in different directions. 
  Chernozukov et al. (2019) and Koike (2019) 
  improved the bound on the growth rate of $p$ and proved validity of the CLT, respectively  for the classes of sets 
  $\A^{max}$ and 
   $\A^{re}$,
  allowing
  $\log p = o(n^{1/5})$. Their  proofs are based on a   randomized version of the Lindeberg's method.  Kuchibolta et al.  (2020) 
   used  techniques from CLT on  Banach Spaces to
  establish  the CLT 
  over $\ell^\infty$-balls in ${\mathcal R}^p$
  (which is a proper subset of $\A^{re}$)
  and also obtained non-uniform large deviation 
  bounds,  allowing $\log p = o(n^{1/4})$.
    Under different sets 
  of additional 
  structural conditions on the $X_i$s (e.g., symmetry and/or existence of a common additive factor along all components), the papers
  by 
  Chernozukov et al. (2019) and Koike (2019) also 
  extend  the 
  CLT over the respective classes of sets,  allowing 
   $\log p = o(n^{1/3})$.
  In all these papers, a key assumption is that the component-wise
  variances of the $X_i$s remain bounded away from zero, which is 
  a critical condition for  anti-concentration of the approximating
  Gaussian measures in ${\mathcal R}^p$.  In the case where this condition fails, 
  it is possible to use the decay of the component-wise variances 
  to reduce the effective dimension of the problem to a lower 
  dimension. Indeed, Lopes et al (2020) consider a similar 
  dimension reduction setting   assuming 
   a polynomial rate of decay of the (ordered)  component-wise variances and, among other interesting results,  establish a uniform CLT over $\A^{max}$ allowing
   $\log p = o(n^{-\delta+ 1/2})$ for any $\delta>0$. 
   Since the maximum is not attained by the low-variance components,
   the assumed decay
   condition on the variances allowed  Lopes et al (2020) to reduce 
   the effective dimension of the problem to a logarithmic 
   scale 
   and 
   apply 
   the classical CLT results in low dimensions growing at  a
   fractional polynomial rate  with  the sample size 
   (cf. 
   Bentkus (2003)).
  In this paper, we do not require   any   such 
  dimension reduction conditions and 
  show that, under some general regularity conditions, 
   \begin{equation}
  \rho_{n,\A^{re}}\raw 0, ~\mbox{provided}~ \log p = o(n^{1/2}) ~\mbox{as}~ \nti. 
   \label{R-1} 
  \end{equation}
  Since $ 
  \A^{max}\subset \A^{re}$, it follows that the best growth rate of $p$ for valid  
  GA over both
  $\A^{max}$
  and 
  $\A^{re}$ is higher than the rate conjectured in [CCK].
  It also improves upon all of the existing results listed above
  allowing a better growth rate of $p$ and matches the rate bound
  in Lopes et al (2020) over a larger class of sets without the 
  variance decay condition.

  Once \eqref{R-1} is established, the CLT result also  raises the
  natural question: 
  {\it When does the GA over 
  $\A^{max}$ or $\A^{re}$ fail in high dimensions?
 }
  In Theorem 2.3, we  show
  by means of an example involving Rademacher random variables 
  that 
  \begin{equation}
  \rho_{n,\A^{re}} \not \raw 0,~~\mbox{as}~ \nti ~\mbox{if}~  
  \limsup_{\nti} ~n^{-(\delta+1/2)} \log p >0,
  \label{R-2} 
  \end{equation}
  for some $\delta>0$. Thus,
  if $\log p$ grows slightly faster than $n^{1/2}$ even along a subsequence, the GA over $\A^{max}$ or $\A^{re}$ in ${\mathcal R}^p$ fails. 
 As a result, the best possible growth rate of $p$
 for a valid  GA over the class $\A^{re}$ is 
 $\log p = o(n^{1/2}) ~\mbox{as}~ \nti$. 
 It can be shown that the  conclusions
 of Theorem 2.3 remain unchanged if we 
 consider $\A=\A^{max}$. Therefore, our main results settle 
 the conjecture of [CCK] by providing
 a definite answer to the critical 
 %on both the upper bound and the lower bound on the 
 growth rate of $p$ for a valid GA.

 The proof of the main results here follows a very different approach compared 
 to [CCK] and other  related recent work in that we make use of 
 the classical Fourier transformation based methods, albeit indirectly. The key tool is
 a set of  non-uniform Berry-Esseen type bounds of Dasgupta (1989, 1992) 
 in the one dimensional CLT for sums of  independent random variables
 which, in turn,  heavily make use of Fourier transformation based 
 arguments 
 %Nagaev () and 
 (cf. Ghosh and Dasgupta (1978)). To derive 
 the GA to $P(T_n\in A)$ over  $A\in \A^{re}$, we begin with the standard factorization of the probability $P(T_n\in A)$ for rectangles $A$ under the (assumed)
 independence of the components of $X_i$ and do a careful analysis of the product
 of the $p$ factors that define $P(T_n\in A)$. Note that  each of the $p$  factors
 must be very close to unity in order to yield a nontrivial value of $P(T_n\in A)$.
 As a result, one must control the 
 errors in the component-wise   normal approximations as well as 
 the probabilities  of the complementary events when the difference from 
 unity is small. We accomplish this by  regrouping the  the endpoints 
 of the $p$ intervals of $A\in \A^{re}$ into a specific  partition of the real line (depending on $n$) 
 and applying a set of  suitable error estimates over  each range. We also needed to  
 make  intricate choices of several tuning parameters to ensure that 
 the final error estimates are close to the actual order of the GA error,
 yielding both a tight upper bound and a lower bound on the growth rate of $p$. 
 See Section \ref{sec:pf} for more details.

%On the central limit theorem in Rp when p → ∞. Prob. Th. %Rel. Fields, 73, 571-583,
%1986

The rest of the paper is organized as follows. We state the assumptions and the main results in Section \ref{sec:main}. 
Proofs of all the results  are presented in Section \ref{sec:pf}.

\section{Main Results}
\label{sec:main}
\setcounter{equation}{0} 
Before moving to the main results, we state 
the assumptions needed to prove the main results.
%\subsection{Assumption}
Suppose $X_1,\dots,X_n$ are independent random vectors in $\mathcal{R}^p$. Let $X_{ij}$ be the $j$th component of $X_i$. Define $s_n^2=\sum_{i=1}^{n}Var(X_{i1})$. Consider the following assumptions:

\begin{enumerate}[label=(A.\arabic*)]
\item $X_{i1},\dots, X_{ip}$ are 
independent and identically distributed 
(iid) for each $i\in \{1,\dots, n\}$.
\item %All order moments of $X_{i1}$ exist for all $i\geq 1$. Additionally,
$\mathbf{E}X_{i1}^{2m-1}=0$ for all $m\geq 1$ and $i\in \{1,\dots,n\}$.
\item %$\liminf_{n\rightarrow \infty}s_n^2>0$.
$0<\inf_{n \geq 1}n^{-1}s_n^2\leq \sup_{n\geq 1}n^{-1}s_n^2<\infty$.
\item $\sup_{n\geq 1}\Big[n^{m-1}\sum_{i=1}^{n}\mathbf{E}\Big(\frac{X_{i1}}{s_n}\Big)^{2m}\Big]\leq \dfrac{l^{-m}(2m)!}{m!}$ for all integer $m\geq 1$, for some $l\in (1,2]$.
\end{enumerate}

Let us discuss the assumptions briefly. 
Under assumption (A.1),  $X_1,\ldots, X_n$ can be non-identically 
distributed (e.g., with a different component-wise 
variance $\sigma^2_i$) but for each fixed $i$, all   $p$ components of $X_i$ 
must have the same distribution. 
The independence of $X_1,\dots,X_n$ and the iid nature among the components of each $X_i$ ensure that $T_1,\dots,T_p$ are iid where $T_j\equiv T_{nj}=s_n^{-1}\sum_{i=1}^{n}X_{ij}$, $j \in \{1,\dots,p\}$. This is essential to keep 
our proof of the main results  simpler. Also note that under (A.1),  $T_n$
in \eqref{clt-1} equals $(T_{n1},\ldots, T_{np})'$. 
Next, to 
gain some insight into assumptions (A.2)-(A.4), consider the 
case when $X_{i}$'s are iid, i.e. when $X_{ij}$'s are all iid. Note that in this case assumption (A.2) is satisfied if $X_{11}$ has a symmetric distribution around $0$ and all moments of $X_{11}$ exist. (A.3) implies and implied by the fact that $X_{11}$ is non-degenerate and has a finite variance. Assumption (A.4) implies that $\mathbf{E}e^{cX^2_{11}}<\infty$ for some $c>0$, which in turn implies that $X_{11}$ has an entire characteristic function. 

Note that in our setup, the distributions of $X_{11},\dots,X_{n1}$ are enough to specify the distribution of $T=s_n^{-1}\sum_{i=1}^{n}X_i$. An immediate example of the sequence $\{X_{11},\dots,X_{n1}\}$ for which all the assumptions are satisfied is when $X_{i1}$'s are iid Rademacher random variables, i.e. when $X_{i1}=1$ or $-1$ each with probability $1/2$. In this case, assumption (A.4) holds with 
$l=2$. For other examples of $\{X_{11},\dots,X_{n1}\}$ which satisfy the above conditions, see section 5 in Dasgupta (1992).  

We are now ready to state the first result. 

\begin{theo}\label{theo:1}
Let $X_1,\dots,X_n$ be independent random vectors in $\mathcal{R}^p$ such that the assumptions (A.1)-(A.4) hold. If $\log p= o(n^{1/2})$ then $$\rho_{n,\mathcal{A}^{re}}\rightarrow 0 \;\;\; \text{as}\;\;n\rightarrow \infty.$$
\end{theo}

Theorem \ref{theo:1} shows that under (A.1)-(A.4),
the GA of  \eqref{clt-2} holds with $\A=\A^{re}$
for $p$ growing at the rate $\exp(o(\sqrt{n}))$
with the sample size $n$. In particular,  \eqref{clt-2} holds 
with $\A=\A^{max}$ beyond the 
range $\log p\gg n^{1/3}$  hypothesized by
[CCK].
%for the validity of the GA for $\A^{max}$,
%and proves the more general version of the 
%conjecture of [CCK]. 
Now a natural question is:  \textit{Does  there exist an asymptotic upper bound on $\rho_{n,\mathcal{A}^{re}}$ even when $\log p$ is exactly  of order $\sqrt{n}$?}  Next theorem gives an answer to this question. 

\begin{theo}\label{theo:2}
Let $X_1,\dots,X_n$ be independent random vectors in $\mathcal{R}^p$ such that the conditions (A.1)-(A.4) hold. Recall that the constant $l\in (1,2]$ is defined in the condition (A.4). Then there exists a positive constant $c\leq (1-l^{-1})^3$ such that whenever $\log p = \epsilon n^{1/2}$ with $0<\epsilon < c$, $$\limsup_{n \rightarrow \infty}\rho_{n,\mathcal{A}^{re}}< \epsilon.$$
\end{theo}

Theorem \ref{theo:2} is a refinement of Theorem \ref{theo:1}
and shows that the uniform error of GA, namely, 
$\rho_{n,\mathcal{A}^{re}}$, decreases linearly with the 
multiplier $\epsilon$ in the rate bound $p \leq \exp(\epsilon \sqrt{n})$ for a nontrivial 
set of $\epsilon$ depending on the tail parameter $l$ 
of the distributions of $X_{11}, \ldots, X_{n1}$. In particular,
when all $X_{ij}$'s are iid with a subGaussian tail, the bound in Theorem \ref{theo:2} holds for all $\epsilon \in (0,c)$ with some $c\leq 1/8$. 
The next result shows that for a valid GA
over the class of sets $\A^{re}$, the $o(\sqrt{n})$ upper 
bound on $\log p$ can not be significantly improved upon.

\begin{theo}\label{theo:3}
Let $X_{ij}$'s be iid Rademacher variables, i.e. $X_{ij}=1$ or $-1$ each with probability $1/2$ and be independent across $i\in \{1\dots,n\}$ and $j\in \{1,\dots, p\}$. If %$\log p > 2\sqrt{n\log n}$
$\limsup_{\nti} n^{-(\delta+1/2)}\log p >0$ for some $\delta>0$, then $$\rho_{n,\mathcal{A}^{re}}\nrightarrow 0\;\; \text{as}\;\; n \rightarrow \infty.$$
\end{theo}

From the proof, it also follows that under the conditions of Theorem \ref{theo:3}, $$\rho_{n,\mathcal{A}^{max}}\nrightarrow 0\;\; \text{as}\;\; n \rightarrow \infty.$$ As a result, the best possible growth rate of $\log p$ for a valid  GA over the 
smaller class of sets $\mathcal{A}^{max}$ considered in [CCK] 
is also $o(\sqrt{n})$.

\section{Proofs of the Results}
\label{sec:pf}
\setcounter{equation}{0} 
Suppose, $\Phi(\cdot)$ and $\phi(\cdot)$ respectively denote the cdf and pdf of the standard normal random variable. Define $N_i=(N_{i1},\dots,N_{ip})^\prime$, $i\in \{1,\dots,n\}$ where $N_{ij}$'s are iid $N(0,1)$ random variables for all $i\in \{1,\dots,n\}$ and $j\in \{1,\dots,p\}$. For any vector $\bm{t}
=(t_1,\dots,t_p)\in \mathcal{R}$, let $t_{(j)}$ and $t^{(j)}$ respectively denote the $j$th element after sorting the components of $\bm{t}$ 
 in increasing order and in decreasing order.
 (We use boldface font only for $\bm{t}$ to avoid some 
 notational conflict later on. All other vectors are 
 denoted using regular font). 
%Also, 
%for vectors $x=(x_1,\ldots,x_p)'$ and   $y=(y_1,\ldots,y_p)'$ in $\bbr^p$,  we write 
%$x\leq y$ if $x_i\leq y_i$ for all $i$.
For any random variable $H$, $P\big(H\leq x\big)$ is assumed to be $1$ if $x=\infty$.  We will need to use some lemmas which are stated are stated and proved next. Proofs of the theorems
are given in Section 3.2 below.

\subsection{Auxiliary Lemmas}
\begin{lem}\label{lem:1}
Let $\{Z_{ni}:1\leq i\leq n\}$ be a triangular array of random variables which are independent within each row with $\sigma_n^2= n^{-1}\sum_{i=1}^{n}\mathbf{E}Z_{ni}^2$ and satisfies
\begin{enumerate}
  \item %All order moments of $X_{i1}$ exist for all $i\geq 1$. Additionally,
$\mathbf{E}Z_{ni}^{2m+1}=0$ for all $m\geq 1$ and $i\in \{1,\dots,n\}$.
\item $0<\inf_{n \geq 1}\sigma_n^2\leq \sup_{n\geq 1}\sigma_n^2<\infty$.
\item $\sup_{n\geq 1}\Big[n^{-1}\sum_{i=1}^{n}\mathbf{E}\Big(\frac{Z_{ni}}{\sigma_n}\Big)^{2m}\Big]\leq \dfrac{q^{-m}(2m)!}{m!}$ for all $m\geq 1$, for some $q\in (1,2]$.
\end{enumerate}
Then we have for some constant $b_1>0$,
$$\Big{|}\mathbf{P}\Big(\sigma_n^{-1}\sum_{i=1}^{n}Z_{ni}\leq t\Big)-\Phi(t)\Big{|}\leq b_1\exp\big(-t^2(1-q^{-1})\big);\;\;\; -\infty<t<\infty.$$
\end{lem}

\noindent
{\bf Proof of Lemma \ref{lem:1}:} 
This lemma is stated as Theorem 1 in Dasgupta (1992).

\noindent
\begin{lem}\label{lem:2}
Under the assumptions of Lemma \ref{lem:1}, 
$$\Big{|}\mathbf{P}\Big(\sigma_n^{-1}\sum_{i=1}^{n}Z_{ni}\leq t\Big)-\Phi(t)\Big{|}\leq b_2r_n\exp\big(-t^2/2\big)
~~\mbox{for all}~ |t|<M_n,
$$
where $M_n =O(n^{1/4})$, $r_n=\max\{n^{-1}M_n^3,n^{-1/2}\}$    and $b_2>0$ is a constant independent of $n$ and $t$.
\end{lem}
\noindent
{\bf Proof of Lemma \ref{lem:2}:}  Let us split $|t|<M_n$ into two parts: when $|t|\leq 1$ and when $1<|t|<M_n$. If $|t|\leq 1$ then by Berry-Esseen theorem (cf. Bhattacharya and Rao (1986)) and assumption (3) of Lemma \ref{lem:1}, we have
\begin{align}\label{eqn:1}
\Big{|}\mathbf{P}\Big(\sigma_n^{-1}\sum_{i=1}^{n}Z_{ni}\leq t\Big)-\Phi(t)\Big{|}&\leq (2.75)\Big[n^{-1}\sum_{i=1}^{n}\mathbf{E}\Big(\frac{|Z_{ni}|}{\sigma_n}\Big)^{3}\Big]n^{-1/2}\nonumber\\
&\leq (2.75)\Big[1+n^{-1}\sum_{i=1}^{n}\mathbf{E}\Big(\frac{Z_{ni}}{\sigma_n}\Big)^{4}\Big]n^{-1/2}\nonumber\\
&\leq \Big[(2.75)\Big(1+\frac{q^{-2}(4!)}{2!}\Big)e^{1/2}\Big]\Big(e^{-t^2/2}n^{-1/2}\Big)\nonumber\\
& \leq b_2r_n\exp(-t^2/2)
\end{align}
Now consider the region $1<|t|<M_n$. Here we are going to use Theorem 2.2 of Dasgupta (1989). Note that under the conditions (2) \& (3) of Lemma \ref{lem:1}, using monotone convergence theorem we have
\begin{align}\label{eqn:2}
&\mathbf{E}\big(Z_{ni}^2\exp(|Z_{ni}|)\big)\leq \mathbf{E}\Big(Z_{ni}^2\big(\exp(|Z_{ni}|+\exp(-|Z_{ni}|)\big)\Big)
= 2\sum_{m=0}^{\infty}\frac{\mathbf{E}|Z_{ni}|^{2m+2}}{(2m)!}\nonumber\\
\Rightarrow\;\;\; &  \sup_{n\geq 1}\bigg[n^{-1}\sum_{i=1}^{n}\mathbf{E}\big(Z_{ni}^2\exp(|Z_{ni}|)\big)\bigg] \leq 2\sum_{m=0}^{\infty}\frac{\sup_{n\geq 1}\big[n^{-1}\sum_{i=1}^{n}\mathbf{E}|Z_{ni}|^{2m+2}\big]}{(2m)!}\nonumber\\
\Rightarrow\;\;\; & \sup_{n\geq 1}\bigg[n^{-1}\sum_{i=1}^{n}\mathbf{E}\big(Z_{ni}^2\exp(|Z_{ni}|)\big)\bigg]\leq 2\sum_{m=0}^{\infty} \dfrac{q^{-(m+1)}(2m+2)! M^{m+1}}{(m+1)!(2m)!}<\infty,
\end{align}
where $M=\sup_{n\geq 1}\sigma_n^2$. Hence in view of (1.3) of Dasgupta (1989), we can consider $g(x)=e^{|x|}$ in applying Theorem 2.2 of Dasgupta (1989). As a consequence, we have for any $1<|t|<M_n$,
\begin{align}\label{eqn:dasgupta89}
\Big{|}\mathbf{P}\Big(\sigma_n^{-1}\sum_{i=1}^{n}Z_{ni}\leq t\Big)-\Phi(t)\Big{|}\leq& b    \exp(-t^2/2)|t|^{-1}|\exp(kn^{-1}t^4)-1|\nonumber\\
&+ b \exp\big(-t^2/2+k n^{-1}t^4\big)n^{-1/2}\nonumber\\
&+ \sum_{i=1}^{n}\mathbf{P}(|Z_{ni}|>r\sqrt{n}\sigma_n|t|)\nonumber\\
=& J_{1n}+J_{2n}+J_{3n}\;\;\;\; \text{(say)},
\end{align}
for some $0<r<1/2$ and constants $b,k>0$ depend on only $r$. Since $n^{-1}t^4=O(1)$, $|\exp(kn^{-1}t^4)-1||t|^{-1}\leq k_1 n^{-1}|t|^3\leq k_1 r_n$ for some constant $k_1>0$, and $\exp(kn^{-1}t^4)=O(1)$. Hence $J_{1n}+J_{2n}\leq b k_2 r_n \exp(-t^2/2)$ for some constant $k_2>0$. Again by Markov's inequality and using (\ref{eqn:2}) we have 
\begin{align*}
\sum_{i=1}^{n}\mathbf{P}(|Z_{ni}|>r\sqrt{n}\sigma_n|t|) &\;\leq r^{-2}\sigma_n^{-2} \bigg[n^{-1}\sum_{i=1}^{n}\mathbf{E}\big(Z_{ni}^2\exp(|Z_{ni}|)\big)\bigg]\exp(-r\sqrt{n}\sigma_n|t|)\\
&\;\leq k_3 r_n \exp(-t^2/2)
\end{align*}
whenever $1<|t|<M_n$, for some constant $k_3>0$. Therefore when $1<|t|<M_n$, taking $b_2=(bk_2 +k_3)$ we have
\begin{align}\label{eqn:3}
\Big{|}\mathbf{P}\Big(\sigma_n^{-1}\sum_{i=1}^{n}Z_{ni}\leq t\Big)-\Phi(t)\Big{|}\leq b_2 r_n \exp(-t^2/2).
\end{align}
Now combining (\ref{eqn:1}) and (\ref{eqn:3}), the proof of Lemma \ref{lem:2} is complete.\\

\begin{lem}\label{lem:3}
Let $Y_1,\dots,Y_n$ be a sequence of mean zero independent random vectors in $\mathcal{R}^p$ with $Y_i = (Y_{i1},\dots Y_{ip})$, $i \in \{1,\dots, n\}$ and let  $\{Y_{i1},\dots,Y_{ip}\}$ be
iid for each $i \in \{1,\dots, n\}$ with $d_n^2= n^{-1}\sum_{i=1}^{n}\mathbf{E}Y_{i1}^2<\infty$. Define, $l_1(x)=\max\Big{\{}\mathbf{P}\Big(d_n^{-1}\sum_{i=1}^{n}\big(-Y_{i1}\big) \leq x \Big), \Phi(x)\Big{\}}$,  $d_1(x)=\Big{|}\mathbf{P}\Big(d_n^{-1}\sum_{i=1}^{n}\big(-Y_{i1}\big) \leq x \Big)-\Phi(x)\Big{|}$, $l_2(x)=\max\Big{\{}\mathbf{P}\Big(d_n^{-1}$ $\sum_{i=1}^{n}Y_{i1} \leq x \Big), \Phi(x)\Big{\}}$ and  $d_2(x)=\Big{|}\mathbf{P}\Big(d_n^{-1}\sum_{i=1}^{n}Y_{i1} \leq x \Big)-\Phi(x)\Big{|}$. Then we have
\begin{align*}
&\Big{|}\mathbf{P}\Big(d_n^{-1}\sum_{i=1}^{n}Y_{i} \in \prod_{j=1}^{p}\big{\{}[a_j,b_j]\cap \mathcal{R}\big{\}}\Big)-\mathbf{P}\Big(n^{-1/2}\sum_{i=1}^{n}N_{i} \in \prod_{j=1}^{p}\big{\{}[a_j,b_j]\cap\mathcal{R}\big{\}}\Big)\Big{|}\\
&\;\leq L_1(\bm{a}) +L_2(\bm{b}),
\end{align*}
where $\bm{a}=(a_1,\dots,a_p)^\prime$, $\bm{b}=(b_1,\dots,b_p)^\prime$, $$L_1(\bm{a})=\bigg[\sum_{k=1}^{p}\Big(\prod_{j\neq k}l_1\big(-a^{(j)}\big)\Big)d_1\big(-a^{(k)}\big)\bigg],\;\;\;\; L_2(\bm{b})=\bigg[\sum_{k=1}^{p}\Big(\prod_{j\neq k}l_2\big(b_{(j)}\big)\Big)d_2\big(b_{(k)}\big)\bigg].$$ 
\end{lem}

\noindent
{\bf Proof of Lemma \ref{lem:3}:}
Note that $d_n^{-1}\sum_{i=1}^{n}Y_{i}=(W_1,\ldots, W_p)'$
where $W_j=d_n^{-1}\sum_{i=1}^{n}$ $Y_{ij}$, $j\in \{1,\dots,p\}$, Then, using the conditions on $Y_{ij}$'s, it is easy to check that $\{W_1,\dots,W_p\}$ are identically distributed.
%(S_1,\dots, S_p)^\prime$ where $S_j=d_n^{-1}\sum_{i=1}^{n}Y_{ij}$, $j\in \{1,\dots,p\}$, and $W_1,\dots,W_p$ are independent. Again $\{W_{i1},\dots,W_{ip}\}$ are identically distributed for all $i \in \{1,\dots, n\}$. Hence $\{W_1,\dots,W_p\}$ are identically distributed. 
Therefore the components of $d_n^{-1}\sum_{i=1}^{n}Y_{i}$, %i.e. $\{W_1,\dots,W_p\}$,
are iid. Similarly, since $N_{ij}\sim N(0,1)$ are iid, the 
$p$-variables $\Big(n^{-1/2}\sum_{i=1}^{n}N_{i1}\Big),\dots,$ $\Big(n^{-1/2}\sum_{i=1}^{n}N_{ip}\Big)$ are also  iid. Hence we have
\begin{align*}
&\Big{|}\mathbf{P}\Big(d_n^{-1}\sum_{i=1}^{n}Y_{i} \in \prod_{j=1}^{p}\big{\{}[a_j,b_j]\cap \mathcal{R}\big{\}}\Big)-\mathbf{P}\Big(n^{-1/2}\sum_{i=1}^{n}N_{i} \in \prod_{j=1}^{p}\big{\{}[a_j,b_j]\cap\mathcal{R}\big{\}}\Big)\Big{|}\\
& = \bigg|\prod_{j=1}^{p}P\Big(W_j \in [a_j,b_j]\cap \mathcal{R}\Big)-\prod_{j=1}^{p}P\Big(n^{-1/2}\sum_{i=1}^{n}N_{ij} \in [a_j,b_j]\cap \mathcal{R}\Big)\bigg|\\
& \leq \bigg[\sum_{k=1}^{p}\Big(\prod_{j\neq k}\Big(\min\Big{\{}l_1\big(-a_j\big),l_2\big(b_j\big)\Big{\}}\Big)\Big)\Big[d_1\big(-a_k\big)+d_2\big(b_k\big)\Big]\bigg]\\
& \leq \sum_{k=1}^{p}\Big(\prod_{j\neq k}l_1\big(-a_j\big)\Big)\Big[d_1\big(-a_k\big)\Big]+ \sum_{k=1}^{p}\Big(\prod_{j\neq k}l_2\big(b_j\big)\Big)\Big[d_2\big(b_k\big)\Big]\\
& = \sum_{k=1}^{p}\Big(\prod_{j\neq k}l_1\big(-a^{(j)}\big)\Big)\Big[d_1\big(-a^{(k)}\big)\Big]+ \sum_{k=1}^{p}\Big(\prod_{j\neq k}l_2\big(b_{(j)}\big)\Big)\Big[d_2\big(b_{(k)}\big)\Big]
\end{align*}
The last equality is due to the following fact:\\
If $(G_1,H_1),\dots,(G_p,H_p)$ are iid random vectors in $\mathcal{R}^2$, then for any $t_1\dots,t_p \in \mathcal{R}$,
\begin{align*}
&\sum_{k=1}^{p}\bigg[\Big(\prod_{j \neq k}\Big(\max\Big{\{}P\Big(G_j \leq t_j\Big), P\Big(H_j\leq t_j\Big)\Big{\}}\Big)\Big)\Big|P\Big(G_k\leq t_k\Big)-P\Big(H_k\leq t_k\Big)\Big|\bigg]\\
&=
\sum_{k=1}^{p}\bigg[\Big(\prod_{j \neq k}\Big(\max\Big{\{}P\Big(G_1 \leq t_{(j)}\Big), P\Big(H_1\leq t_{(j)}\Big)\Big{\}}\Big)\Big)\Big|P\Big(G_1\leq t_{(k)}\Big)-P\Big(H_1\leq t_{(k)}\Big)\Big|\bigg],
\end{align*}
where $\{t_{(1)}, t_{(2)},\dots, t_{(p)}\}$ are obtained after sorting $\{t_1,\dots,t_p\}$ in increasing order.
Therefore we are done.

\begin{lem}\label{lem:4}
For any $t>0$, $\dfrac{1-\Phi(t)}{\phi(t)}\geq \dfrac{2}{\sqrt{t^2+4}+t}$.
\end{lem}

\noindent
{\bf Proof of Lemma \ref{lem:4}:}
This inequality is proved in Birnbaum (1942).\\

\begin{lem}\label{lem:5}
For any positive integer $m$, $$\sqrt{2\pi}\;m^{m+1/2}e^{-m}\leq m! \leq  m^{m+1/2}e^{-m+1}.$$
\end{lem}
This is the well-known Stirling's formula. See for example Robbins (1955).

\subsection{Proofs of the main results}
{\bf Proof of Theorem \ref{theo:1}:} Suppose $T_n=s_n^{-1}\sum_{i=1}^{n}X_i$ and $S_n=n^{-1/2}\sum_{i=1}^{n}N_i$. Let $T=(T_{n1},\dots,T_{np})^\prime$ and $S_n= (S_{n1},\dots,S_{np})^\prime$. Clearly $T_{nj}$'s are iid and $S_{nj}$'s are iid for $j\in \{1,\dots,p\}$. We can use Lemma \ref{lem:3} with $Y_i=X_i$ for $i\in \{1,\dots,n\}$, to obtain
\begin{align*}
&\Big{|}\mathbf{P}\Big(d_n^{-1}\sum_{i=1}^{n}X_{i} \in \prod_{j=1}^{p}\big{\{}[a_j,b_j]\cap \mathcal{R}\big{\}}\Big)-\mathbf{P}\Big(n^{-1/2}\sum_{i=1}^{n}N_{i} \in \prod_{j=1}^{p}\big{\{}[a_j,b_j]\cap\mathcal{R}\big{\}}\Big)\Big{|}\\
\leq & L_1(\bm{a}) +L_2(\bm{b}),
\end{align*}
where $L_1(\bm{a})$ and $L_2(\bm{b})$ are as defined in Lemma \ref{lem:3}. Since all the assumptions are also satisfied if we replace $\{X_1,\dots,X_n\}$ by $\{-X_1,\dots,-X_n\}$, it is enough to show
\begin{align}\label{eqn:4}
\sup_{t_1\leq t_2\leq \dots \leq t_p} L((t_1,\dots,t_p)^\prime)= \sup_{t_1\leq t_2\leq \dots \leq t_p}\bigg[\sum_{j=1}^{p}\Big(\prod_{j\neq k}l(t_j)\Big)d(t_k)\bigg]=o(1),\;\;\;\; \text{as}\; n \rightarrow \infty.
\end{align}
Here, $l(x)=\max\Big{\{}\mathbf{P}\Big(T_{n1} \leq x \Big), \mathbf{P}\Big(S_{n1} \leq x \Big)\Big{\}} \;\;\; \text{and}\;\;\; d(x)=\Big{|}\mathbf{P}\Big(T_{n1} \leq x \Big)-\mathbf{P}\Big(S_{n1} \leq x\Big)\Big{|}$.
Note that we are done if we can show $L(\bm{t}) =L((t_1,\dots,t_p)^\prime) \leq A_n$ for sufficiently large $n$, where $A_n$ does not depend on $\bm{t}$, and $A_n=o(1)$ as $n\rightarrow \infty$. 

Since $\log p =o(n^{1/2})$, there exists a sequence of positive numbers $a_n$ increasing to $\infty$ such that $\log p =O\Big(a_n^{-3}n^{1/2}\Big)$. Without loss of generality assume $a_n^3=o(n^{1/2})$. Now fix $\bm{t}=(t_1,\dots,t_p)^\prime$ in $\mathcal{R}^p$ such that $t_1\leq t_2 \leq \dots\leq t_p$. Then there exist integers $l_1,l_2,l_3$, depending on $n$, such that $0\leq l_1, l_2,l_3 \leq p$ and
\begin{align}\label{eqn:5}
&t_1\leq t_2 \leq \dots \leq t_{l_1} < -a_n^{-1}n^{1/4}\nonumber\\
-a_n^{-1}n^{1/4} \leq\; & t_{l_1+1}\leq t_{l_1+2}\leq \dots \leq\; t_{l_2} < 1\nonumber\\
1 \leq\; & t_{l_2+1}\leq  t_{l_2+2}\leq \dots \leq t_{l_3} \leq a_n^{-1} n^{1/4}\nonumber\\
a_n^{-1}n^{1/4} <\; & t_{l_3+1}\leq  t_{l_3+2}\leq \dots \leq t_p 
\end{align}
Since $a_n^3 = o(n^{1/2})$ and $a_n \rightarrow \infty$ as $n\rightarrow \infty$, hence due to Lemma \ref{lem:2} and Lemma \ref{lem:4}, we have for sufficiently large $n$,
\begin{align}\label{eqn:6}
l(x) \leq\; & I\Big(x>a_n^{-1}n^{1/4}\Big) + \bigg[1-\dfrac{2\phi(1)}{\sqrt{5}+1}+b_2e^{-1/2}n^{-1/4}\bigg]I\Big(x < 1\Big)\nonumber\\
& +\bigg[1-\dfrac{2\phi(x)}{\sqrt{x^2+4}+x}+b_2e^{-x^2/2}n^{-1/4}\bigg]I\Big(x \in\Big[1,a_n^{-1}n^{1/4}\Big]\Big)\nonumber\\
\leq\; & I\Big(x>a_n^{-1}n^{1/4}\Big) + \bigg[1-\dfrac{\phi(1)}{\sqrt{5}+1}\bigg]I\Big(x < 1\Big)\nonumber\\
& +\bigg[1-(4\pi)^{-1/2}a_n n^{-1/4}e^{-x^2/2}\bigg]I\Big(x \in\Big[1,a_n^{-1}n^{1/4}\Big]\Big),
\end{align}
for any $x\in \mathcal{R}$. $I(\cdot)$ is the indicator function. Again due to Lemma \ref{lem:1} and Lemma \ref{lem:2} we have for any $x\in \mathcal{R}$,
\begin{align}\label{eqn:7}
d(x) \leq \Big[b_1 e^{-x^2(1-l^{-1})}\Big]I\Big(|x|>a_n^{-1}n^{1/4}\Big) + \Big[b_2e^{-x^2/2}n^{-1/4}\Big]I\Big(|x|\leq a_n^{-1}n^{1/4}\Big)
\end{align}
Therefore from equations (\ref{eqn:4})-(\ref{eqn:7}), we have
\begin{align}\label{eqn:8}
L(\bm{t})\leq I_1(\bm{t}) +I_2(\bm{t}) +I_3(\bm{t}) + I_4(\bm{t}),
\end{align}
where 
\begin{align*}
I_1(\bm{t}) = &\bigg(\Big[1-\dfrac{\phi(1)}{\sqrt{5}+1}\Big]^{l_2-1}\bigg)*\bigg(\prod_{j=l_2+1}^{l_3}\Big[1-(4\pi)^{-1/2}a_n n^{-1/4}e^{-t_j^2/2}\Big]\bigg)\\
&*\bigg(\sum_{k=1}^{l_1}b_1e^{-t_k^2(1-l^{-1})}\bigg)*I\Big(l_1\geq 1\Big),
\end{align*}
\begin{align*}
I_2(\bm{t})=&\bigg(\Big[1-\dfrac{\phi(1)}{\sqrt{5}+1}\Big]^{l_2-1}\bigg)*\bigg(\prod_{j=l_2+1}^{l_3}\Big[1-(4\pi)^{-1/2}a_n n^{-1/4}e^{-t_j^2/2}\Big]\bigg)\\
&*\bigg(\sum_{k=l_1+1}^{l_2}b_2n^{-1/4}e^{-t_k^2/2}\bigg)*I\Big((l_2-l_1) \geq 1\Big),
\end{align*}
\begin{align*}
I_3(\bm{t})=&\bigg(\Big[1-\dfrac{\phi(1)}{\sqrt{5}+1}\Big]^{l_2}\bigg)*I\Big((l_3-l_2) \geq 1\Big)\\
&*\bigg(\sum_{k=l_2+1}^{l_3}b_2n^{-1/4}e^{-t_k^2/2}\Big(\prod_{ \substack{j=l_2+1\\ j \neq k}}^{l_3}\Big[1-(4\pi)^{-1/2}a_n n^{-1/4}e^{-t_j^2/2}\Big]\Big)\bigg),
\end{align*}
\begin{align*}
I_4(\bm{t})=&\bigg(\Big[1-\dfrac{\phi(1)}{\sqrt{5}+1}\Big]^{l_2}\bigg)*\bigg(\prod_{j=l_2+1}^{l_3}\Big[1-(4\pi)^{-1/2}a_n n^{-1/4}e^{-t_j^2/2}\Big]\bigg)\\
&*\bigg(\sum_{k=l_3+1}^{p}b_1e^{-t_k^2(1-l^{-1})}\bigg)*I\Big((p-l_3) \geq 1\Big).
\end{align*}
\emph{Bound on $I_1(\bm{t})+I_4(\bm{t})$}:
Note that since $\log p =O\big(a_n^{-3}n^{1/2}\big)$, by looking into (\ref{eqn:5}) we have for some $0<M<\infty$,
\begin{align}\label{eqn:9}
I_1(\bm{t})+I_4(\bm{t}) &\leq \bigg(\sum_{k=1}^{l_1}b_1e^{-t_k^2(1-l^{-1})}\bigg)I\Big(l_1\geq 1\Big)+\bigg(\sum_{k=l_3+1}^{p}b_1e^{-t_k^2(1-l^{-1})}\bigg)I\Big((p-l_3) \geq 1\Big)\nonumber\\
& \leq (l_1+p-l_3)\Big(b_1e^{-a_n^{-2}n^{1/2}(1-l^{-1})}\Big)\nonumber\\
& \leq p\Big(b_1e^{-a_n^{-2}n^{1/2}(1-l^{-1})}\Big)\nonumber\\
&\leq b_1e^{Ma_n^{-3}n^{1/2}-a_n^{-2}n^{1/2}
(1-l^{-1})}\Big)\nonumber\\
&\leq b_1e^{-a_n^{-3}n^{1/2}\big(a_n(1-l^{-1})-M\big)}\nonumber\\
 & = A_{1n}\;\;\; \text{(say)}
\end{align}
\emph{Bound on $I_2(\bm{t})$}:
Let $d^{-1}=\Big[1-\dfrac{\phi(1)}{\sqrt{5}+1}\Big]$. Then
\begin{align}
I_2(\bm{t}) &\leq \bigg(\Big[1-\dfrac{\phi(1)}{\sqrt{5}+1}\Big]^{l_2-1}\bigg)\bigg(\sum_{k=l_1+1}^{l_2}b_2n^{-1/4}e^{-t_k^2/2}\bigg)I\Big((l_2-l_1) \geq 1\Big)\nonumber\\
& \leq b_2l_2d^{-(l_2-1)}n^{-1/4}\nonumber\\
& \leq b_2 d n^{-1/4} \sup_{x>0}\big(xd^{-x}\big). \nonumber
\end{align}
Now, $\sup_{x>0}(xc^{-x}) = (\log c)^{-1}c^{-(\log c)^{-1}}$ for any $c>1$. Therefore we have 
\begin{align}\label{eqn:10}
I_2(\bm{t}) \leq \Big(b_2 d (\log d)^{-1}d^{-(\log d)^{-1}}\Big) n^{-1/4} = A_{2n}\;\;\; \text{(say)}
\end{align}
\emph{Bound on $I_3(\bm{t})$}: Note that if $(l_3-l_2) = 0$ then $I_3(\bm{t})=0$ and there is nothing more to do. Hence assume $(l
_3-l_2)\geq 1$. Then we have
\begin{align} \label{eqn:11}
I_3(\bm{t}) \leq \bigg(\sum_{k=l_2+1}^{l_3}b_2n^{-1/4}e^{-t_k^2/2}\Big(\prod_{ \substack{j=l_2+1\\ j \neq k}}^{l_3}\Big[1-(4\pi)^{-1/2}a_n n^{-1/4}e^{-t_j^2/2}\Big]\Big)\bigg) =I_{31}(\bm{t})\;\;\; \text{(say)}.
\end{align}
We are going to check the monotonicity of $I_{31}(\bm{t})$ with respect to $t_{l_2+1},\dots,t_{l_3}$. Note that 
\begin{align*}
\dfrac{\partial I_{31}(\bm{t})}{\partial t_l }=& \bigg[b_2t_le^{-t_l^2/2}n^{-1/4}\prod_{\substack{j=l_2+1\\ j\neq l}}^{l_3}\Big[1-(4\pi)^{-1/2}a_n n^{-1/4}e^{-t_j^2/2}\Big]\bigg]\times\\
&\bigg[\sum_{\substack{k=l_2+1\\k\neq l}}^{l_3}\bigg(\Big[1-(4\pi)^{-1/2}a_n n^{-1/4}e^{-t_k^2/2}\Big]^{-1}(4\pi)^{-1/2}a_n n^{-1/4}e^{-t_k^2/2}\bigg)-1\bigg]
\end{align*}
Hence for any $l=l_2+1,\dots,l_3$, $\dfrac{\partial I_{31}(\bm{t})}{\partial t_l }\gtreqless 0$ if and only if 
\begin{align}\label{eqn:12}
\sum_{\substack{j=l_2+1 \\ j\neq l}}^{l_3}\dfrac{z_j}{1-z_j}\gtreqless 1,
\end{align}
where $z_j=(4\pi)^{-1/2}a_n n^{-1/4}e^{-t_j^2/2}$ for $j\in \{l_2+1,\dots,l_3\}$. Note that since $1 \leq t_{l_2+1}\leq \dots \leq t_{l_3}$, $1> z_{l_2+1} \geq \dots \geq z_{l_3} >0$ for sufficiently large $n$. Hence for sufficiently large $n$, $$\dfrac{z_{l_2+1}}{1-z_{l_2+1}} \geq \dots \geq \dfrac{z_{l_3}}{1-z_{l_3}},$$ due to the fact that $z/(1-z)$ is increasing for $z\in (0,1)$. Therefore from (\ref{eqn:12}) we can say that $I_{31}(\bm{t})$ is non-increasing in $\{t_{l_2+1},\dots,t_m\}$ and non-decreasing in $\{t_{m+1},\dots,t_{l_3}\}$ where $(m-l_2)$ is a non-negative integer not more than $(l_3-l_2)$. Again note that $1\leq t_{l_2+1}\leq \dots \leq t_m \leq t_{m+1} \leq \dots t_{l_3}\leq a_n^{-1}n^{1/4}$. Hence from (\ref{eqn:11}) we have
\begin{align*}
I_3(\bm{t})\leq I_{31}((\bm{t}^{(1)\prime},\bm{t}^{(2)\prime})^\prime)
\end{align*}
where $\bm{t}^{(1)}$ is an $(m-l_2)\times 1$ vector with each component being $1$ and $\bm{t}^{(2)}$ is an $(l_3-m)\times 1$ vector with each component being $a_n^{-1}n^{1/4}$. Therefore using the fact that $\log p =O(a_n^{-3}n^{1/2})$ we can say that there exists $M\in(0,\infty)$ such that for sufficiently large $n$,
\begin{align}\label{eqn:13}
I_3(\bm{t}) \leq\; & (m-l_2)\Big[1-(4\pi)^{-1/2}a_n n^{-1/4}e^{-1/2}\Big]^{m-l_2-1}\big(b_2n^{-1/4}e^{-1/2}\big)\nonumber\\
&+ (l_3-m)\big(b_2n^{-1/4}e^{-(2^{-1}a_n^{-2}n^{1/2})}\big)\nonumber\\
\leq \;& 2\Big[\exp\Big(\log (m-l_2)-(m-l_2)\big((4\pi)^{-1/2}a_n n^{-1/4}e^{-1/2}\big)\Big)\Big]\big(b_2n^{-1/4}e^{-1/2}\big)\nonumber\\
& + p\big(b_2n^{-1/4}e^{-(2^{-1}a_n^{-2}n^{1/2})}\big)\nonumber\\
\leq\;&  2\Big[\exp\Big(\sup_{x>0}\big[\log x-x\big((4\pi)^{-1/2}a_n n^{-1/4}e^{-1/2}\big)\big]\Big)\Big]\big(b_2n^{-1/4}e^{-1/2}\big)\nonumber\\
&+ \Big(b_2n^{-1/4}e^{-\big((a_n^{-3}n^{1/2})(a_n/2-M)\big)}\Big)\nonumber\\
 \leq\; & \big((4\pi)^{1/2}a_n^{-1} n^{1/4}e^{1/2}\big)\big(b_2n^{-1/4}e^{-1/2}\big) + \Big(b_2n^{-1/4}e^{-\big((a_n^{-3}n^{1/2})(a_n/2-M)\big)}\Big)\nonumber\\
=\;&A_{3n}\;\;\; \text{(say)}
\end{align}
Combining (\ref{eqn:9}), (\ref{eqn:10}) and (\ref{eqn:13}), we have for sufficiently large $n$,
$$I(\bm{t})\leq A_{1n} + A_{2n} +A_{3n}=A_n\;\;\;\text{(say)}.$$ Note that $A_n$ does not depend on the choice of $\bm{t}$ and also $A_n\rightarrow 0$ as $n\rightarrow \infty$. Therefore, the proof of Theorem \ref{theo:1} is now complete.\\

\noindent
{\bf Proof of Theorem \ref{theo:2}:}  We are going to follow the same steps as in the proof of Theorem \ref{theo:1}. Note that we are done if in (\ref{eqn:4}) we can show $L_n(\bm{t})\leq \epsilon/2$ for sufficiently large $n$, irrespective of the choice of $\bm{t}$.

Now take $c=\min\Big\{[8b_2\sqrt{2\pi}(\sqrt{2}+1)]^{-3},(1-l^{-1})^3\Big\}$ where the constant $b_2$ is as defined in Lemma \ref{lem:2} but with $\{Z_{ni}:1\leq i\leq n\}=\{X_1,\dots,X_n\}$ and the constant $l$ is defined in the condition (A.4). Fix $\bm{t}=(t_1,\dots,t_p)^\prime$ in $\mathcal{R}^p$ such that $t_1\leq t_2 \leq \dots\leq t_p$. Then there exist integers $l_4,l_5,l_6$, depending on $n$, such that $0\leq l_4, l_5,l_6 \leq p$ and
\begin{align}\label{eqn:14}
&t_1\leq t_2 \leq \dots \leq t_{l_4} < -\epsilon^{1/3} n^{1/4}\nonumber\\
-\epsilon^{1/3}n^{1/4} \leq\; & t_{l_4+1}\leq t_{l_4+1}\leq \dots \leq\; t_{l_5} < 1\nonumber\\
1 \leq\; & t_{l_5+1}\leq  t_{l_5+1}\leq \dots \leq t_{l_6} \leq \epsilon^{1/3} n^{1/4}\nonumber\\
\epsilon^{1/3}n^{1/4} <\; & t_{l_6+1}\leq  t_{l_6+1}\leq \dots \leq t_p 
\end{align}

Now use the same definitions of $l(x)$ and $u(x)$, as in the proof of Theorem \ref{theo:1}. Then due to Lemma \ref{lem:1} \& \ref{lem:2} and the fact that $\epsilon \leq c < 1$, we have
\begin{align}\label{eqn:15}
l(x) \leq\; & I\Big(x>\epsilon^{1/3}n^{1/4}\Big) + \bigg[1-\dfrac{2\phi(1)}{\sqrt{5}+1}+\epsilon b_2 e^{-1/2}n^{-1/4}\bigg]I\Big(x < 1\Big)\nonumber\\
& +\bigg[1-\dfrac{2\phi(x)}{\sqrt{x^2+4}+x}+\epsilon b_2 e^{-x^2/2}n^{-1/4}\bigg]I\Big(x \in\Big[1,\epsilon^{1/3}n^{1/4}\Big]\Big)\nonumber\\
\leq\; & I\Big(x>\epsilon^{1/3}n^{1/4}\Big) + \bigg[1-\dfrac{\phi(1)}{\sqrt{5}+1}\bigg]I\Big(x < 1\Big)\nonumber\\
& +\bigg[1-\dfrac{n^{-1/4}e^{-x^2/2}}{\sqrt{2\pi}\big(\sqrt{\epsilon^{2/3}+4n^{-1/2}}+\epsilon^{1/3}\big)} \bigg]I\Big(x \in\Big[1,\epsilon^{1/3}n^{1/4}\Big]\Big),
\end{align}
and
\begin{align}\label{eqn:16}
u(x) \leq \Big[b_1 e^{-x^2(1-l^{-1})}\Big]I\Big(|x|>\epsilon^{1/3}n^{1/4}\Big) + \Big[b_2\epsilon n^{-1/4}e^{-x^2/2}n^{-1/4}\Big]I\Big(|x|\leq \epsilon^{1/3}n^{1/4}\Big)
\end{align}
for any $x\in \mathcal{R}$, for sufficiently large $n$. Therefore from equations (\ref{eqn:4})-(\ref{eqn:7}) we have for sufficiently large $n$,
\begin{align}\label{eqn:17}
L(\bm{t})\leq J_1(\bm{t}) +J_2(\bm{t}) +J_3(\bm{t}) + J_4(\bm{t}),
\end{align}
where 
\begin{align*}
J_1(\bm{t}) =\; & \bigg(\Big[1-\dfrac{\phi(1)}{\sqrt{5}+1}\Big]^{l_5-1}\bigg)*\bigg(\prod_{j=l_5+1}^{l_6}\Big[1-\dfrac{n^{-1/4}e^{-t_j^2/2}}{\sqrt{2\pi}\big(\sqrt{\epsilon^{2/3}+4n^{-1/2}}+\epsilon^{1/3}\big)} \Big]\bigg)\\
&*\bigg(\sum_{k=1}^{l_4}b_1e^{-t_k^2(1-l^{-1})}\bigg)*I\Big(l_4\geq 1\Big),
\end{align*}
\begin{align*}
J_2(\bm{t})=\;&\bigg(\Big[1-\dfrac{\phi(1)}{\sqrt{5}+1}\Big]^{l_5-1}\bigg)*\bigg(\prod_{j=l_5+1}^{l_6}\Big[1-\dfrac{n^{-1/4}e^{-t_j^2/2}}{\sqrt{2\pi}\big(\sqrt{\epsilon^{2/3}+4n^{-1/2}}+\epsilon^{1/3}\big)} \Big]\bigg)\\
&*\bigg(\sum_{k=l_3+1}^{l_4}\epsilon b_2n^{-1/4}e^{-t_k^2/2}\bigg)I\Big((l_2-l_1) \geq 1\Big),
\end{align*}
\begin{align*}
J_3(\bm{t})=\;&\bigg(\Big[1-\dfrac{\phi(1)}{\sqrt{5}+1}\Big]^{l_5}\bigg)*I\Big((l_6-l_5) \geq 1\Big)\\
&*\bigg(\sum_{k=l_5+1}^{l_6}\epsilon b_2n^{-1/4}e^{-t_k^2/2}\Big(\prod_{ \substack{j=l_5+1\\ j \neq k}}^{l_6}\Big[1-\dfrac{n^{-1/4}e^{-t_j^2/2}}{\sqrt{2\pi}\big(\sqrt{\epsilon^{2/3}+4n^{-1/2}}+\epsilon^{1/3}\big)} \Big]\Big)\bigg),
\end{align*}
\begin{align*}
J_4(\bm{t})=\;&\bigg(\Big[1-\dfrac{\phi(1)}{\sqrt{5}+1}\Big]^{l_5}\bigg)*\bigg(\prod_{j=l_5+1}^{l_6}\Big[1-\dfrac{n^{-1/4}e^{-t_j^2/2}}{\sqrt{2\pi}\big(\sqrt{\epsilon^{2/3}+4n^{-1/2}}+\epsilon^{1/3}\big)} \Big]\bigg)\\
&*\bigg(\sum_{k=l_6+1}^{p}b_1e^{-t_k^2(1-l^{-1})}\bigg)I\Big((p-l_6) \geq 1\Big).
\end{align*}

\emph{Bound on $J_1(\bm{t})+J_4(\bm{t})$}:
Since $\log p =\epsilon n^{1/2}$, from (\ref{eqn:17}) we have
\begin{align}\label{eqn:18}
J_1(\bm{t})+J_4(\bm{t}) &\leq \bigg(\sum_{k=1}^{l_4}b_1e^{-t_k^2(1-l^{-1})}\bigg)I\Big(l_4\geq 1\Big)+\bigg(\sum_{k=l_6+1}^{p}b_1e^{-t_k^2(1-l^{-1})}\bigg)I\Big((p-l_6) \geq 1\Big)\nonumber\\
& \leq p\Big(b_1e^{-\epsilon^{2/3}n^{1/2}(1-l^{-1})}\Big)\nonumber\\
&= b_1\exp{\big(\epsilon n^{1/2}-\epsilon^{2/3}n^{1/2}
(1-l^{-1})\big)}\nonumber\\
&\leq b_1\exp{\big(-\epsilon^{2/3}n^{1/2}\big((1-l^{-1})-\epsilon^{1/3}\big)}\nonumber\\
 & < \epsilon/12,
\end{align}
for large enough $n$, since $\epsilon < c \leq (1-l^{-1})^3$.

\emph{Bound on $J_2(\bm{t})$}:
Noting that $d^{-1}=\Big[1-\dfrac{\phi(1)}{\sqrt{5}+1}\Big]$ and $\epsilon<1$, we have for sufficiently large $n$,
\begin{align}\label{eqn:19}
J_2(\bm{t})\leq I_2(\bm{t}) &\leq \bigg(\Big[1-\dfrac{\phi(1)}{\sqrt{5}+1}\Big]^{l_2-1}\bigg)\bigg(\sum_{k=l_1+1}^{l_2}b_2n^{-1/4}e^{-t_k^2/2}\bigg)I\Big((l_2-l_1) \geq 1\Big)\nonumber\\
&\leq  \Big(b_2d^{-1} (\log d)^{-1}d^{-(\log d)^{-1}}\Big) n^{-1/4} < \epsilon/12.
\end{align}

\emph{Bound on $J_3(\bm{t})$}: Write $z_n = n^{-1/4}e^{-1/2}\Big[\sqrt{2\pi}\big(\sqrt{\epsilon^{2/3}+4n^{-1/2}}+\epsilon^{1/3}\big)\Big]^{-1}$. Through the same line of arguments as in bounding $I_3(\bm{t})$, we have
for some non-negative integer $q \in [l_5,l_6]$ and for sufficiently large $n$,
\begin{align}\label{eqn:20}
J_3(\bm{t}) \leq\;& (q-l_5)\Big[1-z_n \Big]^{q-l_5-1}\big(\epsilon b_2n^{-1/4}e^{-1/2}\big) \nonumber\\
&+ (l_6-q)\big(\epsilon b_2n^{-1/4}e^{-(2^{-1}\epsilon^{2/3}n^{1/2})}\big)\nonumber\\
\leq\;& 2\Big[\exp\Big(\log (q-l_5)-(q-l_5)z_n\Big)\Big]\big(\epsilon b_2n^{-1/4}e^{-1/2}\big)\nonumber\\
&+ p\big(\epsilon b_2n^{-1/4}e^{-(2^{-1}\epsilon^{2/3}n^{1/2})}\big)\nonumber\\
\leq\;&  2\Big[\exp\Big(\sup_{x>0}\big[\log x-xz_n\big]\Big)\Big]\big(\epsilon b_2n^{-1/4}e^{-1/2}\big)\nonumber\\
&+ \Big(\epsilon b_2n^{-1/4}\exp{\big(-\epsilon^{2/3}n^{1/2}\big(1/2-\epsilon^{1/3}\big)}\Big)\nonumber\\
\leq\;& 2z_n^{-1}\big(\epsilon b_2n^{-1/4}e^{-1/2}\big) + \Big(\epsilon b_2n^{-1/4}\exp{\big(-\epsilon^{2/3}n^{1/2}\big(1/2-\epsilon^{1/3}\big)}\Big)\nonumber\\
\leq\;& \Big[\sqrt{2\pi}\big(\sqrt{\epsilon^{2/3}+4n^{-1/2}}+\epsilon^{1/3}\big)\Big]2\epsilon b_2\nonumber\\
&+ \Big(\epsilon b_2n^{-1/4}\exp{\big(-\epsilon^{2/3}n^{1/2}\big(1/2-\epsilon^{1/3}\big)}\Big)\nonumber\\
 <\;& \epsilon/4 +\epsilon/12,
\end{align}
since $\epsilon^{1/3} < c^{1/3}\leq [8b_2\sqrt{2\pi}(\sqrt{2}+1)]^{-1}$. Now combining (\ref{eqn:14})-(\ref{eqn:20}), the proof of Theorem \ref{theo:2} is complete.\\

\noindent
{\bf Proof of Theorem \ref{theo:3}:} Suppose
$\limsup_{\nti} n^{-(\delta+1/2)} \log p >0 $
for some $\delta>0$. Then, there exists a subsequence
$\{n_k\}$ such that $\log p_{n_k}  > n_k^{\delta+1/2}$ for all 
$k\in \bbn$. % and for some $\delta\in (0,1)$. 
We will  consider  two cases  
depending on the values of $p$:\\
(I) $p_{n_k}   \geq  \sqrt{\dfrac{\pi}{2}}\Big(\sqrt{n_k+4}+\sqrt{n_k}\Big)e^{n_k/2}
$ for infinitely many $k\in \bbn$;\\
(II) $ n_k^{\delta+1/2}
< \log p_{n_k} \leq 
 \log\Big[\sqrt{\dfrac{\pi}{2}}\Big(\sqrt{n_k+4}+\sqrt{n_k}\Big)e^{n_k/2}\Big]
$ for all but finitely many $k\in \bbn$. \\

%Note that 
 %$\delta \in (0, 1]$ under Case (II). 
Next, strictly for the sake of notational simplicity, without 
loss of generality, we shall suppose that the respective inequalities under
Cases (I) and (II) hold for all $n\in \bbn$. 
(Otherwise,
one needs to 
extract a further subsequence $\{n_{k_i}\}$ of $\{n_k\}$
for Case (I) and rewrite all the steps below for Case(I) 
with $n$ replaced by $n_{k_i}$, and do
similarly for Case (II)).\\

\noindent
\underline {Case (I):}
Suppose that $p\geq \sqrt{\dfrac{\pi}{2}}\Big(\sqrt{n+4}+\sqrt{n}\Big)e^{n/2}$
for all $n\in \bbn$.  We will show that $\rho_{n,\mathcal{A}^{re}}\nrightarrow 0$, as $n\rightarrow \infty$.  Note that here $s_n^2=n$. Consider the set $A=(-\infty,\sqrt{n}]^p$. Then clearly $A\in \mathcal{A}^{re}$ and
\begin{align}\label{eqn:21}
P\Big(s_n^{-1}\sum_{i=1}^{n}X_i \in A\Big) =  \Big[P\Big(\sum_{i=1}^{n}X_{i1} \leq n\Big)\Big]^p=1.  
\end{align}
However using Lemma \ref{lem:4} we have
\begin{align}\label{eqn:22}
P\Big(n^{-1/2}\sum_{i=1}^{n}N_i \in A\Big) &=  \Big[P\Big(n^{-1/2}\sum_{i=1}^{n}N_{i1} \leq \sqrt{n}\Big)\Big]^p  \nonumber\\
& = \Big[\Phi\Big(\sqrt{n}\Big)\Big]^p \nonumber\\
& \leq \Big[1-\dfrac{2\phi(\sqrt{n})}{\sqrt{n+4}+\sqrt{n}}\Big]^p \nonumber\\
& = \Big[1-\dfrac{2 e^{-n/2}}{\sqrt{2\pi}\big(\sqrt{n+4}+\sqrt{n}\big)}\Big]^p.
\end{align}

Therefore if $p\geq \sqrt{\dfrac{\pi}{2}}\Big(\sqrt{n+4}+\sqrt{n}\Big)e^{n/2}$ then $P\Big(s_n^{-1}\sum_{i=1}^{n}N_i \in A\Big) \leq 2/e$ for large enough $n$ and hence from (\ref{eqn:21}) it is clear that $\rho_{n,\mathcal{A}^{re}}\nrightarrow 0$, as $n\rightarrow \infty$, when $p\geq \sqrt{\dfrac{\pi}{2}}\Big(\sqrt{n+4}+\sqrt{n}\Big)e^{n/2}$.\\

\noindent
\underline{Case (II):}
Now suppose that 
$n^{1/2+\delta}\leq \log p\leq \log\Big[\sqrt{\dfrac{\pi}{2}}\Big(\sqrt{n+4}+\sqrt{n}\Big)e^{n/2}\Big]$ for all $n\in \bbn$. Note that  here
$\delta\in (0,1)$.
We will show that $\rho_{n,\mathcal{A}^{re}}\nrightarrow 0$, as $n\rightarrow \infty$. Consider the set $B= \Big(-\infty, n^{1/4}f(n)\Big]^p$. Here $\{f(n)\}_{n\geq 1}$ is a sequence such that $n^{3/4}f(n)$ is an even integer and $n^{\delta/4}\leq f(n)\leq \dfrac{n^{1/4}}{1+\eta}$ with some constant $0<\eta < 1$. $\eta$ is going to be specified later. For rest of the proof, we will only consider $n$ to be even.
(The proof is similar for the odd integer subsequence, with 
 the lower limit of the summation in \eqref{eqn:23}
 changed to $[n_1+1]/2$).
Then writing $n_1=n^{3/4}f(n)$ we have
\begin{align}\label{eqn:23}
L_{1n}=P\Big(s_n^{-1}\sum_{i=1}^{n}X_i \in B\Big) &=  \Big[P\Big(\sum_{i=1}^{n}X_{i1} \leq n^{3/4}f(n)\Big)\Big]^p\nonumber\\
& = \bigg[1- \sum_{\substack{k=n_1+2 \\ k \;\text{is an even integer}}}^{n} \binom n {\frac{(n-k)}{2}} 2^{-n}\bigg]^p\nonumber\\
& = \bigg[1- \sum_{k=\frac{n_1}{2} + 1}^{\frac{n}{2}} \binom n {\frac{n}{2}-k} 2^{-n}\bigg]^p\nonumber\\
&\geq \bigg[1-\Big(\frac{n-n_1}{2}\Big)\binom n {\frac{n-n_1}{2}} 2^{-n}\bigg]^p\nonumber\\
& = \bigg[1-g(n)\bigg]^p\;\;\;\; \text{(say)},
\end{align}
where the last inequality follows due to the fact that $$\binom n 1 \leq \binom n {2} \leq \dots \leq \binom n {l+1},
\;\;\;\;\; \text{for any positive integer}\;\; l < \frac{n-1}{2}.$$
Now applying Lemma \ref{lem:5} we have
\begin{align*}
g(n) \leq \frac{(n-n_1)e(2\pi)^{-1}n^{-1/2}}{\big(1-\frac{n_1}{n}\big)^{\frac{n-n_1+1}{2}}\big(1+\frac{n_1}{n}\big)^{\frac{n+n_1+1}{2}}},
\end{align*}
where writing $t_1=\Big(\frac{n-n_1+1}{2}\Big)\log \big(1-\frac{n_1}{n}\big)$ and $t_2=\Big(\frac{n+n_1+1}{2}\Big)\log \big(1+\frac{n_1}{n}\big)$ we have
\begin{align*}
(t_1+t_2)&=\sum_{k=1}^{\infty}\frac{n_1^{2k}}{n^{2k-1}}\Big(\frac{1}{2k-1}-\frac{1}{2k}\Big) + \frac{1}{2}\log \Big(1-\frac{n_1^2}{n^2}\Big)\\
& \geq \frac{n_1^2}{2n} + \frac{n_1^4}{14n^3} + \frac{1}{2}\log \Big(1-\frac{1}{(1+\eta)^2}\Big),
\end{align*}
for sufficiently large $n$. Therefore we have for sufficiently large $n$,
\begin{align}\label{eqn:24}
L_{1n}&\geq \bigg[1-\Big[e (2\pi)^{-1}\Big(1-\frac{1}{(1+\eta)^2}\Big)^{1/2}\Big]\Big(\frac{n-n_1}{\sqrt{n}}\Big) e^{-\frac{n_1^2}{2n} - \frac{n_1^4}{14n^3}}\bigg]^p\nonumber\\
& = \big[1-E_{1n}\big]^p\;\;\;\; \text{(say)}.
\end{align}
Again note that by Lemma \ref{lem:4},
\begin{align}\label{eqn:25}
L_{2n}=P\Big(n^{-1/2}\sum_{i=1}^{n}N_i \in B\Big) &=  \Big[P\Big(n^{-1/2}\sum_{i=1}^{n}N_{i1} \leq n^{-1/2}n_1\Big)\Big]^p\nonumber\\
& \leq \Bigg[1-\frac{2\phi\Big(\frac{n_1}{\sqrt{n}}\Big)}{\sqrt{\frac{n_1^2}{n}+4}+\frac{n_1}{\sqrt{n}}}\Bigg]^p\nonumber\\
&= \Bigg[1-\frac{\big(\pi/2\big)^{-1/2}e^{-\frac{n_1^2}{2n}}}{\sqrt{\frac{n_1^2}{n}+4}+\frac{n_1}{\sqrt{n}}}\Bigg]^p\nonumber\\
& = \big[1-E_{2n}\big]^p\;\;\;\; \text{(say)}.
\end{align}
Now observe that $
\frac{n^{\delta}}{14}\leq \frac{n_1^4}{14n^3}\leq \frac{n}{(1+\eta)^4}$, due to the assumed condition $n^{\delta/4}\leq f(n)\leq \dfrac{n^{1/4}}{1+\eta}$. Again for the choice of $f(n)$, $\log 3n_1 <  n^{\delta/2}$ for large enough $n$. Therefore 
$$E_{1n}\ll E_{2n}.$$
We shall choose  $f(n)$ suitably depending on the 
growth rate of $p$  to show
that $\rho_{n,\mathcal{A}^{re}}\nrightarrow 0$ as $n\rightarrow \infty$. When $n^{\delta+1/2} < \log p < \Big[\dfrac{3\log n}{4} + \dfrac{n}{2(1+\eta)^2}\Big]$, it easy to choose a sequence $f(n)$ in the given range such that $pE_{2n}\raw \infty$ and $pE_{1n}\raw 0$,
so that $L_{1n}\raw 1$ but $L_{2n}\raw 0$ and 
$\rho_{n,\mathcal{A}^{re}}\nrightarrow 0$ as $n\rightarrow \infty$.  In particular take $\sqrt{n} [f(n)]^2 \approx 
2[\log p - \log n]$.

The case when $\log p$ is of comparable order to $n$, or more precisely when $ \Big[\dfrac{3\log n}{4} + \dfrac{n}{2(1+\eta)^2}\Big] \leq \log p \leq 
 \log\Big[\sqrt{\dfrac{\pi}{2}}\Big(\sqrt{n+4}+\sqrt{n}\Big)e^{n/2}\Big]$, the choice of $f(n)$ is a little tricky. 
We take $f(n)=\frac{n^{1/4}}{1+\eta}$ where $\eta>0$ is to be 
specified below. Then, it is easy to check that 
$p E_{2n}\raw \infty$. 
Next we show that for large enough $n$, $E_{1n}^{-1}\gg \sqrt{\dfrac{\pi}{2}}\Big(\sqrt{n+4}+\sqrt{n}\Big)e^{n/2}$ for some  $\eta>0$. % when $f(n)=\frac{n^{1/4}}{1+\eta}$.
 This is true if $\frac{(1+\eta)n}{2}\leq \frac{n}{2(1+\eta)^2}+\frac{n}{14(1+\eta)^4}$, that is if 
\begin{align}\label{eqn:26}
\bigg[\frac{1}{(1+\eta)^3}+\frac{1}{7(1+\eta)^5}-1\bigg]\geq 0.
\end{align}
Fix an $\eta>0$ such that this condition holds. 
Then it follows that $pE_{1n}\rightarrow 0$, as $n\rightarrow \infty$ as before and hence $\rho_{n,\mathcal{A}^{re}}\nrightarrow 0$ as $n\rightarrow \infty$. 
This completes the proof of Theorem
\ref{theo:3}.

\bibliographystyle{amsplain}

\end{document}